\documentclass{elsart}

\usepackage{amsmath}
\usepackage{amssymb}
\newtheorem{question}[thm]{Question}

\begin{document}

\begin{frontmatter}

\title{Problems related to\\ a de Bruijn - Erd\H{o}s theorem}

\author{Xiaomin Chen}

\address{Department of Computer Science,\\
 Rutgers University, Piscataway, NJ 08854-8019, USA}

\author{Va\v sek Chv\'atal\thanksref{VC}}

\thanks[VC]{This research was undertaken, in part, thanks to funding
  from the Canada Research Chairs Program and from the Natural
  Sciences and Engineering Research Council of Canada.}

\address{Canada Research Chair in Combinatorial Optimization,\\  
Department of Computer Science and Software Engineering,\\
Concordia University, Montr\' eal, Qu\' ebec H3G 1M8, Canada}

\bigskip

\author{\em In memory of Leo Khachiyan}

\begin{abstract}
De Bruijn and Erd\H os proved that every noncollinear set of $n$
points in the plane determines at least $n$ distinct lines.  We
suggest a possible generalization of this theorem in the framework of
metric spaces and provide partial results on related extremal
combinatorial problems.
\end{abstract}

\begin{keyword}
combinatorial geometry \sep metric space \sep metric betweenness \sep
extremal combinatorial problem

\end{keyword}

\end{frontmatter}

\section{Lines in metric spaces}

Two distinct theorems are referred to as ``the de Bruijn - Erd\H os
theorem''. One of them \cite{DE51} concerns the chromatic number of
infinite graphs; the other \cite{DE48} is our starting point: {\em
Every noncollinear set of $n$ points in the plane determines at least
$n$ distinct lines.\/}

This theorem involves neither measurement of distances nor measurement
of angles: the only notion employed here is incidence of points and
lines. Such theorems are a part of {\em ordered geometry\/}
\cite{C61}, which is built around the ternary relation of {\em
betweenness\/}: point $y$ is said to lie between points $x$ and $z$ if
$y$ is an interior point of the line segment with endpoints $x$ and
$z$. It is customary to write $[xyz]$ for the statement that $y$ lies
between $x$ and $z$. In this notation, a {\em line\/} $\overline{uv}$
is defined \ ---\  for any two distinct points $u$ and $v$ \ ---\  as
\begin{equation}\label{def.first}
\{p: [puv]\}\;\cup\; \{u\}\;\cup\; \{p: [upv]\}\;\cup\; \{v\}\;\cup\; \{p: [uvp]\}. 
\end{equation}

In terms of the Euclidean metric $\rho$, we have
\begin{equation}\label{def.btw}
\mbox{ $[abc] \;\Leftrightarrow\; a,b,c$ are three distinct points and
$\rho(a,b)+\rho(b,c) = \rho(a,c)$.}
\end{equation}
For an arbitrary metric space, equivalence (\ref{def.btw}) defines the
ternary relation of {\em metric betweenness} introduced in \cite{Men}
and further studied in \cite{Blu,Bus,Chv}; in turn, (\ref{def.first})
defines the line $\overline{uv}$ for any two distinct points $u$ and
$v$ in the metric space. The resulting family of lines may have
strange properties. For instance, a line can be a proper subset of
another: in the metric space with points $u,v,x,y,z$ and 
\begin{eqnarray}
\hspace{1.5cm} &&
\rho(u,v)=\rho(v,x)=\rho(x,y)=\rho(y,z)=\rho(z,u)=1,\nonumber \\
&& \rho(u,x)=\rho(v,y)=\rho(x,z)=\rho(y,u)=\rho(z,v)=2,\nonumber
\end{eqnarray}
we have
\[
\overline{vy}=\{v,x,y\} \;\;\mbox{ and }\;\;
\overline{xy}=\{v,x,y,z\}. 
\]

Nevertheless, fragments of ordered geometry might translate to the
framework of metric spaces. In particular, we know of no
counterexample to the de Bruijn - Erd\H os theorem in this framework.

\begin{question}\label{open}
True or false? Every finite metric space $(X,\rho)$ where no line
consists of the entire ground set $X$ determines at least $|X|$
distinct lines.
\end{question}

\section{Lines in hypergraphs}

A {\em hypergraph\/} is an ordered pair $(X,H)$ such that $X$
is a set and $H$ is a family of subsets of $X$; elements of $X$
are the {\em vertices\/} of the hypergraph and members of $H$
are its {\em edges\/}. Our definition of lines in a metric space
$(X,\rho)$ depends only on the hypergraph $(X,H(\rho))$ where 
\[
H(\rho)=\{\{a,b,c\}:\,[abc]\}:
\]
the line $\overline{uv}$ equals $\{u,v\}\cup\{p:\{u,v,p\}\in
H(\rho)\}.$ This observation leads us to extend the notion of lines in
metric spaces to a notion of lines in hypergraphs: given an arbitrary
hypergraph $(X,H)$, we define the line $\overline{uv}$ \ ---\ for any
two distinct vertices $u$ and $v$ \ ---\ as $\{u,v\}\cup\{p: \;
\exists T \, (T\in H, \{u, v, p\}\subseteq T) \,\}.$ Now every
metric space $(X,\rho)$ and its associated hypergraph $(X,H(\rho))$
define the same family of lines.

A hypergraph is called {\em $k$-uniform\/} if each of its edges
consists of $k$ vertices. All the hypergraphs $(X,H(\rho))$ are
$3$-uniform, but some $3$-uniform hypergraphs do not arise from any
metric space $(X,\rho)$ as $(X,H(\rho))$: it has been proved
(\cite{Chv,Che}) that the {hypergraph consisting of the seven vertices
$0,1,2,3,4,5,6$ and the seven edges
\[
\{ i \bmod{7},\, (i+1) \bmod{7},\, (i+3) \bmod{7} \}\hspace{1cm} (i=0,1,2,3,4,5,6)
\]
does not arise from any metric space. (This $3$-uniform hypergraph is
known as the {\em Fano plane\/} or the {\em projective plane of order
two\/}.) Restricting the notion of lines to $3$-uniform hypergraphs
would bring about no loss of generality: for every hypergraph $(X,H)$
there is a $3$-uniform hypergraph $(X,H^{(3)})$ such that $(X,H)$ and
$(X,H^{(3)})$ define the same family of lines. Specifically,
\[
H^{(3)} = \{S: |S|=3 \;\mbox{ and }\; \exists T \, (T\in H, S\subseteq T) \,\}.
\]

Let $m(n,k)$ denote the smallest number of lines in a hypergraph on
$n$ vertices where every line consists of at most $k$
vertices. Showing that $m(n,n-1)\ge n$ would show that the answer to
Question \ref{open} is ``true''. However, as we are going to prove,
$m(n,n-1)$ grows slower than every power of $n$.

\begin{lem}\label{construct}
If $n,\ell ,a$ are positive integers such that
$
2\le n-\ell  \le a^{\,\ell} ,
$
then 
\[
m(n,n-1)\le 2^{\,\ell} +\ell a.
\]
\end{lem}

\begin{pf}
Write $P=\{1,2,\ldots ,\ell \}$ and let $A$ be a set of size $a$. By
assumption, there is a set $S$ of strings of length $\ell $ over alphabet
$A$ such that $|S|=n-\ell $ and such that, for each $i$ in $P$, some two
strings in $S$ differ in their $i$-th position. For each choice of $i$
in $P$ and $x$ in $A$, set
\[
E_{ix}=\{i\}\cup\{x_1x_2\ldots x_\ell \in S: x_i=x\}.
\]
Now consider all the lines $\overline{uv}$ in the hypergraph
\[
(P\cup S,\: \{P,S\}\cup\{E_{ix}: i\in P, x\in A\}).
\]
If $u,v\in P$, then $\overline{uv}=P$.  If $u\in P$ and $v\in S$, then
$\overline{uv}=E_{ux}$ with $x$ the $u$-th character in $v$. If
$u,v\in S$, then $\overline{uv}=S\cup P'$ with $P'$ the set of
positions in which $u$ and $v$ agree; $P'$ is a proper (and possibly
empty) subset of $P$. So the hypergraph has $n$ vertices, none of its
lines consists of all $n$ vertices, and there are at most $1+ \ell a +
(2^{\,\ell} -1)$ lines. \hskip1em $\Box$
\end{pf}
 
\begin{thm}
There are positive constants $n_0$ and $c$ such that
\begin{equation}\label{ub}
n\ge n_0 \;\Rightarrow\; m(n,n-1)\le c^{\sqrt{\ln n}}.
\end{equation}
for all $n$.
\end{thm}

\begin{pf}
Let $\alpha,\beta,\gamma,\delta$ be arbitrary constants such that
\[
0<\alpha<1<\beta<\gamma<2<\delta.
\]
There is a positive integer $\ell _0$ such that
\[
\ell \ge \ell _0 \;\;\Rightarrow\;\; 
\alpha \ell <\ell -1,\;  
\beta^{\,\ell} <\gamma^{\,\ell} -1,\; 
\ell \gamma^{\,\ell} <2^{\,\ell} ,\;
2^{\,\ell +1}<\delta^{\,\ell} \!.
\]
We claim that (\ref{ub}) holds as long as
\[
n\ge n_0 \;\Rightarrow\; n-\left\lceil\sqrt\frac{\ln   n}{\ln\beta}\,\right\rceil\ge 2 
\]
and 
\[
\ln n_0\ge \ell _0^2\ln\beta,\;\;\;
\ln c\ge  \frac{\ln\delta}{\alpha\sqrt{\ln \beta}}\, .
\]
To justify this claim, consider an arbitrary $n$ such that $n\ge n_0$
and set
\[
\ell =\left\lceil\sqrt\frac{\ln n}{\ln\beta}\,\right\rceil, \;\;\;
a=\lfloor\gamma^{\,\ell} \rfloor.
\]
Now $\ell \ge \ell _0$, $a>\beta^{\,\ell} $, and so $\ell \ln a>\ell ^2\ln\beta\ge \ln
n$. Lemma \ref{construct} guarantees that
\[
m(n,n-1)\le 2^{\,\ell} +\ell a;
\]
since
\[
\ell <\frac{\ell -1}{\alpha}<\frac{1}{\alpha}\sqrt\frac{\ln n}{\ln\beta}
\]
we have
\[
2^{\,\ell} +\ell a < 2^{\,\ell +1} <\delta^{\,\ell} <c^{\sqrt{\ln n}}.
\]
$\Box$
\end{pf}

We do not know the order of growth of $m(n,n-1)$; our best lower bound
is only logarithmic in $n$. (We follow the convention of letting $\lg$
stand for the logarithm to base $2$.)

\begin{thm}\label{lb1}
$m(n,n-1)\ge \lg n$.
\end{thm}

\begin{pf}
Consider an arbitrary hypergraph with $n$ vertices and $m$ lines where
no line consists of all $n$ vertices. Let us observe that
\begin{eqnarray}
&&\mbox{for every two distinct vertices $u$ and $v$,} \nonumber \\ 
&&\mbox{there is a line which includes $u$ and does not include $v$: } \label{antichain}
\end{eqnarray}
by assumption, some vertex $w$ is not included in line
$\overline{uv}$, and so no edge includes all three vertices $u,v,w$,
and so line $\overline{uw}$ includes $u$ and does not include $v$. 
For each vertex $x$, let $S_x$ denote the set of all lines that
include $x$. Property (\ref{antichain}) guarantees that these $n$ sets
are all distinct, and so $n\le 2^m$. \hskip1em $\Box$
\end{pf}

Actually, property (\ref{antichain}) guarantees that the $n$ sets
$S_x$ form an {\em antichain\/} in the sense that none of them is a
subset of another. This observation allows a negligible improvement of
the bound in Theorem \ref{lb1}: first, the classic result of Sperner
(\cite{Spe}) asserts that an antichain on a ground set of size $m$ has
at most
\[
\binom{m}{\lfloor m/2 \rfloor}
\]
sets; next, by Stirling's formula,
\[
\binom{m}{\lfloor m/2 \rfloor} \sim \frac{2^m}{\sqrt{\pi m/2}}\, ;
\]
finally, if $m=\lg n + \frac{1}{2}\lg\lg n +c$, then
\[
\frac{2^m}{\sqrt{\pi m/2}} \sim 2^c(2/\pi)^{1/2}n.
\]
It follows that for every positive $\varepsilon$ there is an $n_0$
such that
\[
n\ge n_0 \;\Rightarrow\; m(n,n-1)>
\lg n + \frac{1}{2}\lg\lg n + \frac{1}{2}\lg\frac{\pi}{2}- \varepsilon.
\]

Since $m(n,k)$ is a nonincreasing function of $k$, Theorem~\ref{lb1}
guarantees that $m(n,k)\ge \lg n$ whenever $2\le k < n$. For
small values of $k$, this bound can be much improved.

\begin{thm}\label{lb2}
\[
m(n,k)\ge \frac{n(n-1)}{k(k-1)}
\]
whenever $n\ge k\ge 2$. 
\end{thm}

\begin{pf}
Consider an arbitrary hypergraph with $n$ vertices and $m$ lines where
every line consists of at most $k$ vertices. Trivially,
\begin{eqnarray}
&&\mbox{for every two distinct vertices $u$ and $v$,} \nonumber \\
&&\mbox{there is a line which includes both $u$ and $v$. }
  \label{cover}
\end{eqnarray}
Let $P$ denote the set of all pairs $(L,\{u,v\})$ such that $L$ is a
line and $u,v$ are two distinct vertices in $L$. On the one hand,
every line includes at most $k$ points, and so
\[
|P|\le m {k\choose 2}.
\]
On the other hand, property (\ref{cover}) guarantees that
\[
|P|\ge {n\choose 2}.
\]
The lower bound on $m$ follows by comparing the two bounds on
$|P|$. \hskip1em $\Box$
\end{pf}

When the value of $k$ is fixed, the lower bound of Theorem \ref{lb2}
is asymptotically optimal:

\begin{thm}\label{ub2}
\[
\lim_{n\rightarrow\infty}\: m(n,k)\cdot\frac{k(k-1)}{n(n-1)}\,=\,1
\]
whenever $k\ge 2$.
\end{thm}

\begin{pf}
Theorem \ref{lb2} guarantees that
\[
\liminf_{n\rightarrow\infty}\: m(n,k)\cdot\frac{k(k-1)}{n(n-1)}\,\ge\,1\, .
\]
In every $k$-uniform hypergraph $(X,H)$ such that
\begin{equation}\label{packing}
\mbox{every two edges share at most one vertex,}
\end{equation}
each line is either an edge or a set of two vertices that is not a
subset of any edge, and so there are 
\[
|H|\;+\;\left({|X|\choose 2}-|H|{k\choose 2}\right)
\]
lines altogether. In particular, with $f(n,k)$ standing for the
largest number of edges in a $k$-uniform hypergraph with $n$ vertices
and with property (\ref{packing}), we have
\[
m(n,k)\le {n\choose 2}-f(n,k)\left({k\choose 2}-1\right);
\]
Erd{\H{o}}s and Hanani \cite{EH} proved that
\[
\lim_{n\rightarrow\infty}\: f(n,k)\cdot\frac{k(k-1)}{n(n-1)}\,=\,1\, ;
\]
it follows that 
\[
\limsup_{n\rightarrow\infty}\: m(n,k)\cdot\frac{k(k-1)}{n(n-1)}\,\le\,1\, .
\]
$\Box$
\end{pf}

\section{Closure-lines in hypergraphs and metric spaces}

The {\em Sylvester-Gallai theorem\/} \cite{Syl,Erd,C61,BM,EP,PA,Chv}
asserts that every noncollinear finite set $X$ of points in the plane
includes two points such that the line passing through them includes
no other point of $X$. This theorem does not translate to the
framework of metric spaces along the simple lines of our Section~1: in
the five-point example of that section, every line consists of three
or four ponts. Nevertheless, it does translate to the framework of
metric spaces in a circuitous way, which we are about to describe.

Let us call a set $T$ of vertices in a hypergraph {\em affinely
closed\/} if, and only if, every edge that shares at least two
vertices with $T$ is fully contained in $T$. For every set $S$ of
vertices, the intersection of all affinely closed supersets of $S$ is
an affinely closed set, which we will refer to as the {\em affine
closure\/} of $S$ and which we will denote by ${\rm aff}(S)$.  By {\em
closure-lines\/} in the hypergraph, we shall mean all the sets ${\rm
aff}(\{u,v\})$ with $u$ and $v$ two distinct vertices; by
closure-lines in a metric space $(X,\rho)$, we shall mean
closure-lines in its associated hypergraph $(X,H(\rho))$.

When $X$ is a subset of a Euclidean space and $\rho$ is the Euclidean
metric, lines and closure-lines in $(X,\rho)$ coincide: each of them
is the intersection of $X$ and the Euclidean line passing through two
distinct points of $X$. One of us \cite{Chv} conjectured and the other
one \cite{Che} proved that the notion of closure-lines provides a
translation of the Sylvester-Gallai theorem to the framework of metric
spaces:
\begin{quote}
{\em In every finite metric space, some closure-line includes\\
 either all the points of the ground set or only two of them.}
\end{quote}
The same notion falls far short of providing a translation of the de
Bruijn - Erd\H{o}s theorem to the framework of metric spaces:
\begin{thm}
For every integer $n$ greater than $5$, there is a metric space on $n$
points where each closure-line consists of at most $n-2$ points and
there are precisely $7$ distinct closure-lines altogether.
\end{thm}
\begin{pf}
Consider the metric space $(X,\rho)$, where  
$X=\{x_k:1\le k\le n\}$ with 
\begin{eqnarray}
\mbox{\hspace{1.6cm}}& x_1=(1,3),\;\;\; x_2=(2,4),\;\;\;
x_3=(3,1),\;\;\; x_4=(4,2),&
\nonumber \\
& x_k=(k,n+5-k)\;\mbox{whenever}\;5\le k\le n,& \nonumber 
\end{eqnarray}
and
\[
\rho((a_1,a_2),\,(b_1,b_2))\;=\; |a_1-b_1|+|a_2-b_2|.
\]
Since $H(\rho)$ consists of all $\{x_1,x_2,x_k\}$ with $5\le k\le n$,
all $\{x_3,x_4,x_k\}$ with $5\le k\le n$, and all $\{x_i,x_j,x_k\}$
with $5\le i<j<k\le n$, we have 
\begin{eqnarray}
{\rm aff}(\{x_1,x_2\})&=&X-\{x_3,x_4\},\nonumber \\
{\rm aff}(\{x_3,x_4\})&=&X-\{x_1,x_2\},\nonumber \\
{\rm aff}(\{x_i,x_j\})&=&X-\{x_1,x_2,x_3,x_4\}\;\,\mbox{whenever $5\le
  i<j\le n$},\nonumber \\
{\rm aff}(\{x_i,x_j\})&=&\{x_i,x_j\}\;\,\mbox{whenever $1\le i\le 2$
  and $3\le j\le 4$},\nonumber \\
{\rm aff}(\{x_i,x_j\})&=&X-\{x_3,x_4\}\;\,\mbox{whenever $1\le i\le 2$
  and $5\le j\le n$},\nonumber \\
{\rm aff}(\{x_i,x_j\})&=&X-\{x_1,x_2\}\;\,\mbox{whenever $3\le i\le 4$ 
and $5\le j\le n$}.\nonumber 
\end{eqnarray}
$\Box$
\end{pf}

Finally, let $\overline{m}(n,k)$ denote the smallest number of
closure-lines in a hypergraph on $n$ vertices where every closure-line
consists of at most $k$ vertices. Our proof of Theorem \ref{lb2} with
``lines'' replaced by ``closure-lines'' shows that
\begin{equation}\label{lb3}
\overline{m}(n,k)\ge \frac{n(n-1)}{k(k-1)}
\end{equation}
whenever $n\ge k\ge 2$;$\;$ in turn, our proof of Theorem \ref{ub2} with
``lines'' replaced by ``closure-lines'' yields the following conclusion.
\begin{thm}\label{limit}
\[
\lim_{n\rightarrow\infty}\:
\overline{m}(n,k)\cdot\frac{k(k-1)}{n(n-1)}\,=\,1\,
\]
whenever $k\ge 2$.
\end{thm}
The order of growth of $\overline{m}(n,k)$ is given by its lower bound
(\ref{lb3}):
\begin{thm}\label{sandwich}
There is a positive constant $c$ such that
\[
\frac{n(n-1)}{k(k-1)}\;\le \;\overline{m}(n,k)\;\le \;
c\cdot\frac{n(n-1)}{k(k-1)}
\]
whenever $n\ge k\ge 2$. 
\end{thm}
\begin{pf}
For every integer $k$ greater than $1$, Theorem \ref{limit} guarantees
the existence of a constant $c_k$ such that
\begin{equation}\label{small}
\overline{m}(n,k)\le
c_k\cdot\frac{n(n-1)}{k(k-1)}\;\;\mbox{whenever $n\ge k$.}
\end{equation}
With $c$ any constant such that
\[
c\ge 12\;\;\;\mbox{and}\;\;\; c\ge c_k\;\mbox{whenever $2\le k< 12$},
\]
we propose to show that, for every integer $k$ greater than $1$, 
\begin{equation}\label{ub3}
\overline{m}(n,k)\;\le \;
c\cdot\frac{n(n-1)}{k(k-1)}\;\;\mbox{whenever $n> k$.}
\end{equation}
(Trivially, $\overline{m}(n,k)=1$ whenever $2\le n \le k$.)  For this
purpose, consider an arbitrary but fixed integer $k$ greater than
$1$. If $k<12$, then (\ref{ub3}) follows from (\ref{small}); if $k\ge
12$, then we will use induction on $n$ to prove that $\overline{m}(n,k)
\le cn^2/k^2$ whenever $n>k$.

Set
\[
p=2\left\lceil\frac{n+1}{k}\right\rceil\
\]
and note for a future reference that
\[
4\le p<2\left(\frac{n+1}{k}+1\right)\le 
\frac{4n}{k}\: .
\]
Take a set $X$ such that $|X|=n$, take a subset $X_0$ of $X$ such that
$|X_0|=p-1$, and partition $X-X_0$ into pairwise disjoint sets $V_i\;
(1\le i\le p)$ whose sizes are as nearly equal as possible. Since
\[
\frac{k}{4}-1 < \frac{n-(p-1)}{p} \le \frac{k}{2}-1\, ,
\]
we have
\[
2\le \min |V_i| \le \max |V_i|\le \frac{k-1}{2}\: .
\]
In some hypergraph $(X_0, H_0)$, every closure-line consists of at
most $k$ vertices and there are precisely $\overline{m}(p-1,k)$
distinct closure-lines altogether.  A theorem of Behzad, Chartrand,
and Cooper, Jr. \cite{BCC} guarantees that (the chromatic index of the
complete graph $K_{2s}$ is $2s-1$, and so) there is a mapping
\[
\phi:\{S:\, S\subset \{1,2,\ldots ,p\},\, |S|=2\} \rightarrow X_0
\]
with the following property:
\begin{quote}
for every $i$ in $\{1,2,\ldots ,p\}$ and for every $w$ in $X_0$\\
there is precisely one $j$ in $\{1,2,\ldots ,p\}$ such that $\phi(\{i,j\})=w$.
\end{quote}
Set
\begin{eqnarray}
H_1 &=& \{\{u,v,w\}:\; 
\mbox{there are $i$ and $j$ with}\; u\in V_i,\, v\in V_j,\, \phi(\{i,j\})=w\},
\nonumber \\
H_2 &=& \{S:\; |S|=3 \;\mbox{and there is an $i$ with}\; S\subseteq V_i\},\nonumber 
\end{eqnarray}
and $H=H_0\cup H_1\cup H_2$. Since closure-lines in hypergraph $(X,H)$
are
\begin{itemize}
\item all the closure-lines in hypergraph $(X_0,H_0)$,
\item all the sets $V_i\cup V_j\cup \{\phi(\{i,j\})\}$ such that $1\le
  i<j\le p$, and
\item all the sets $V_i$ such that $1\le i\le p$,
\end{itemize}
we have
\[
\overline{m}(n,k) \le \overline{m}(p-1,k) + {p\choose 2} +p\, .
\]
If $p-1>k$, then (as $p-1<n/3$) the induction hypothesis guarantees
that 
\[
\overline{m}(p-1,k) \le c\,\left(\frac{p-1}{k}\right)^2 <
\frac{c}{9}\left(\frac{n}{k}\right)^2\! ;
\]
if $p-1\le k$, then 
\[
\overline{m}(p-1,k) =1 < \frac{c}{9}\left(\frac{n}{k}\right)^2\! ;
\]
finally,
\[
{p\choose 2} +p = {p+1\choose 2} < 10\left(\frac{n}{k}\right)^2\! .
\]
We conclude that 
\[
\overline{m}(n,k) \le \frac{c}{9}\left(\frac{n}{k}\right)^2 +
10\left(\frac{n}{k}\right)^2 \le c\cdot\frac{n^2}{k^2}\: .
\]
$\Box$
\end{pf}

\begin{ack}
We are grateful to Bahman Kalantari for a question which stimulated
our development of Theorem~\ref{lb1}. 
\end{ack}

\end{document}